\documentclass[10pt]{article}
\usepackage[cp1251]{inputenc}
\usepackage[english,russian]{babel}
\usepackage{amsmath,amssymb,amsthm}
\usepackage{chngcntr,color}
\usepackage{ulem}
\definecolor{navyblue}{rgb}{0.0, 0.0, 0.5} 
\usepackage[colorlinks=true, citecolor = navyblue,  linkcolor=navyblue, filecolor=magenta, urlcolor=navyblue]{hyperref}
\usepackage[mathscr]{euscript}
\usepackage{enumitem}
\usepackage{cite}
\usepackage{bbm} 

\makeatletter
\def\namedlabel#1#2{\begingroup
    #2%
    \def\@currentlabel{#2}%
    \phantomsection\label{#1}\endgroup
}
\makeatother

\newtheorem{theorem}{Теорема}
\newtheorem{lemma}{Лемма}[section]

\newtheorem{corollary}{Следствие}

\numberwithin{equation}{section}

\newcommand{\ra}{\rightarrow}
\newcommand{\p}[1]{{\mathbf{P}} \left( \, #1 \, \right) }

\newcommand{\e}{{\mathbf{E}} }

\newcommand{\skk}[1]{\left\{ #1 \right\}}

\renewcommand{\le}{\leqslant}
\renewcommand{\ge}{\geqslant}

\def\v{\varepsilon}

\def\d{\delta}
\def\a{\alpha}

\def\g{\gamma}
\def\l{\lambda}
\def\m{\mu}
\def\({\left(}
\def\){\right)}

\def\({\left(}
\def\){\right)}

\def\[{\left[}
\def\]{\right]}

\def\SL{\left\{}
\def\SP{\right\}}

\def\aa{{\boldsymbol{\alpha}}}
\def\bb{{\boldsymbol{\beta}}}

\def\mm{{\boldsymbol{\mu}}}
\def\0{{\boldsymbol{0}}}

\def\n{\nu}
\def\t{\tau}

\def\z{\zeta}
\def\zz{{\boldsymbol{\zeta}}}
\def\e{\eta}
\def\h{\theta}

\def\mm{{\boldsymbol{\mu}}}
\def\xx{{\boldsymbol{\xi}}}

\def\R{\mathbb R}
\def\Z{\mathbb Z}

\newcommand{\scalar}[2]{\left\langle   #1, #2 \right\rangle }

\begin{document}

\title{Large deviations principle for terminating multi\-di\-men\-sional compo\-und
renewal processes with application to polymer pinning models
}
 \vspace{3cm}
 \author{ A. Logachov, A. Mogulskii, E. Prokopenko  
 }
 \date{}
 \maketitle
 
\textbf{Annotation}. Large deviations principle is obtained for terminating multi\-di\-men\-sional compo\-und
renewal processes. We also obtained the asymptotic of large deviations 
for the case when a Gibbs change of the original probability measure takes place.
The random processes mentioned in the paper are widely used in  polymer pinning models.


\begin{center}
\textbf{Принцип больших уклонений для многомерных  обобщенных процессов восстановления  с приложением
к связыванию полимеров}
\end{center}
 
\begin{center}
А.В. Логачев, А.А. Могульский, Е.И. Прокопенко
\end{center}
 
\textbf{Аннотация}. В работе получен принцип больших уклонений для обрывающихся многомерных обобщенных
процессов восстановления. Также получена асимптотика больших уклонений, для случая,
когда происходит гиббсовская замена исходной вероятностной меры.
Рассмотренный тип случайных процессов  широко используется в моделях
связывания полимеров.

\textbf{Ключевые слова} обобщенный процесс восстановления, принцип больших уклонений,
функционал уклонений,
модель связывания полимеров, гиббсовская замена меры.

\tableofcontents

\section{Введение}
Работа посвящена  изучению предельного поведения вероятностной меры
построенной по многомерному обобщенному процессу  восстановления (ОПВ), который,
вообще говоря, может обрываться (т.е. допускается возможность того, что
время между моментами восстановления может быть равным $\infty$ с положительной вероятностью).
Такие случайные процессы находят свое применение, в частности, в моделях
связывания полимеров (polymer pinning models). Прежде чем привести
обзор известных результатов, нам удобней дать
строгое математическое определение
изучаемого объекта.

Пусть случайный вектор $(\t,\zz,v)$ принимает  значения в пространстве
$$
\overline{\R}_+\times \R^d\times\R,
$$
где $\overline{\R}_+:=\{t\in \R:~t>0\}\cup\{\infty\}$,
так что координата $\t>0$ может принимать значение $\infty$  c вероятностью
 ${\bf P}(\t=\infty)=:p\in [0,1)$. Далее, пусть
 $$
   \{(\t_i,\zz_i,v_i)\}_{i\ge 1}
 $$
 --- последовательность независимых копий вектора $(\t,\zz,v)$.
 Обозначим
 $$
 T_0:=0,~~T_n:=\sum_{i=1}^n\t_i,~~~ \mathbf{Z}_0:=\mathbf{0},~~\mathbf{Z}_n:=\sum_{i=1}^n\zz_i,~~~
  V_0:=0,~~V_n:=\sum_{i=1}^nv_i,
 $$
где $n\in \Z_+$.

По суммам  $(T_n,\mathbf{Z}_n)$ для $t>0$ построим однородный ОПВ  с $d$-мерным фазовым
пространством    $\R^d$:
$$
  {\mathbf{Z}}(t):={\mathbf{Z}}_{\n(t)},~~~\n(t):=\sup\{n\ge 0:~T_n\le t\},
$$
возможно обрывающийся (в случае, когда $p={\bf P}(\t=\infty)>0$).

Рассмотрим  семейство
\begin{equation*}\label{def_prob_Marco}
  {\bf P}_t(\mathbf{Z}(t)\in B):=\frac{{\bf E}\(e^{V_{\n(t)}};~\mathbf{Z}(t)\in B\)}
  {{\bf E}e^{V_{\n(t)}}},~~~~B\in {\cal B}_d,
\end{equation*}
вероятностных распределений в пространстве $\R^d$, где
${\cal B}_d$ --- борелевская $\sigma$-алгебра в $\R^d$.

Нас будет интересовать асимптотическое поведение последовательности
\begin{equation}{\label{n.1}}
\frac{1}{t}\ln {\bf P}_t\(\frac{\mathbf{Z}(t)}{t}\in B\)
\end{equation}
при $t\to \infty$.

Сделаем теперь краткий обзор известных результатов, связанных с асимптотикой последовательности (\ref{n.1}).
Работы, в которых изучается предельное поведение последовательности (\ref{n.1})
можно разделить на две группы.
 Первая носит чисто теоретический характер,
в ней изучается асимптотика, связанная с поведением непосредственно
обобщенного процесса восстановления, т.е. полагают, что
$V_{\n(t)}\equiv 0$ и изучают асимптотику последовательности
\begin{equation}\label{a04.11.1}
\frac{1}{t}\ln {\bf P}\(\frac{\mathbf{Z}(t)}{t}\in B\).
\end{equation}
В этом направлении, для случая $\mathbf{P}(\tau=\infty)=0$, хорошо изучена как грубая асимптотика
последовательности (\ref{a04.11.1}) (принцип больших уклонений (ПБУ))
\cite{MogProk3,MogProk5,Tsi}, так и точная асимптотика
(локальные и инегро-локальные теоремы в области нормальных, умеренно больших и больших уклонений)
 \cite{BorMog2} (для одмерного случая),
 \cite{MogProk1,MogProk2} (для многомерного случая), \cite{MogProkLog1}
(для многомерных арифметических полумарковских ОПВ).
В работах
\cite{MogProk4,MogLog1} изучена, соответственно, грубая и точная асимптотика больших уклонений для конечномерных приращений многомерных ОПВ. Работы \cite{BorMog1,MogLog2,Mog1}
посвящены принципам больших и умеренно больших уклонений для траекторий ОПВ.
Отметим также работу \cite{LefMarZam},
в ней получен ПБУ для мер построенных по ОПВ.
В работе \cite{Bak}  для обрывающихся многомерных ОПВ
в некоторой области фазового пространства установлена интегро-локальная
(локальная) предельная теорема.

Вторая группа работ имеет более прикладной характер, в ней изучается
асимптотическое поведение последовательности (\ref{n.1}) для случая, когда,
вообще говоря, $V_{\n(t)}\neq 0$ \cite{Zamp1,Zamp2}. В этих работах
$\tau$ --- целое число, которое является случайным количеством
мономеров присоединяемых к имеющемуся полимеру в процессе синтеза,
$\zz$ --- числовая характеристика, присоединяемого блока мономеров
(например, количества мономеров того или иного вида в присоединяемом блоке).
В частности, если $\tau$ --- количество мономеров, присоединенных
до выхода полимера на границу разных сред, то вероятность того,
с какой стороны от границы был присоединен блок мономеров
зависит от количества энергии в этих средах.
(см. например \cite[глава 1]{Gant}, \cite[глава 2]{Holl}). Таким образом, из физических соображений возникает
потребность рассматривать не исходную вероятностную меру (\ref{a04.11.1}),
а некоторое ее экспоненциальное преобразование (\ref{n.1}).
Отметим также, что в этих моделях, в частности,  событие $\left\{\tau=\infty\right\}$ может означать невозможность
присоединения нового блока мономеров.

Следует отметить, что в работах из первой группы тоже приходится сталкиваться с экспоненциальным
преобразованием исходной вероятностной меры, но совершенно из других
соображений, связанных с техникой доказательства ПБУ.
Эта техника помогает также успешно решить многие задачи,
поставленные во второй группе работ.

Наиболее близкой по содержанию к нашей работе является
статья \cite{Zamp2}.
В ней рассматривается случайный процесс $\mathbf{Z}(t)$ у которого фазовое пространство является банаховым, в том числе и бесконечномерным, и, в частности, успешно изучается асимптотическое поведение последовательности (\ref{n.1})
для достаточно широкого класса множеств $B \in \mathcal{B}_d$, 
в случае, когда выполнены дополнительные условия:
\begin{itemize}
\item[$(a)$] {\it случайная величина $\t$ принимает целые положительные значения или
$\infty$},
\item[$(b)$] {\it случайная величина $v$ имеет вид
$v=h(\t)$ для некоторой фиксированной неслучайной функции $h=h(t).$}
\end{itemize}

В настоящей работе  выполнено изучение асимптотического поведения последовательности
(\ref{n.1}) в общем случае, т.е. без привлечения условий $(a)$ и $(b)$
(см. следствие \ref{cor_cor1}).
Кроме того, для обрывающихся многомерных ОПВ в условиях, близких к необходимым установлен
принцип больших уклонений в фазовом пространстве (см. следствие \ref{cor_cor2}).
Отметим, что более общий вид функции $v=h(\t)$
дает возможность рассматривать более сложные физические модели синтеза полимеров.
В частности, модели в которых экспоненциальное преобразование меры зависит
от типа присоединенных мономеров. Отказ от целочисленности $\tau$ позволяет, в частности,
рассматривать модели в которых мы можем фиксировать только моменты
времени выхода полимера на границу сред и не знаем какое точно количество
мономеров было присоединено за время между этими моментами.
Вопрос о выполнении ПБУ для ОПВ в произвольном
банаховом пространстве, когда условия $(a)$ и $(b)$ не выполнены, остается
открытым.

Для формулировок и доказательства основных результатов нам потребуются некоторые дополнительные моментные условия.
Везде далее будем предполагать, что выполнено следующее условие
\begin{itemize}
\item[\namedlabel{C}{$[\mathbf{C}^*]$}] Для некоторых $\lambda > 0$ и  $ \varepsilon >0$
$$\mathbf{E}\left(e^{-\lambda \tau+\varepsilon|\zz|+v} ;\, \tau<\infty\right)<\infty.$$
\end{itemize}
Заметим, что в силу того, что $-\l \infty=-\infty$, в условии
\ref{C} можно оставить неравенство
$$
  {\bf E}\(e^{-\l\t+\v|\zz|+ v}\)<\infty.
$$

Поскольку
\begin{equation}\label{16.12.1}
 \ln {\bf P}_t\(\frac{\mathbf{Z}(t)}{t}\in B\)
 = \ln {\bf E}\(e^{V_{\n(t)}};~\frac{\mathbf{Z}_{\n(t)}}{t}\in B\)
 -\ln {\bf E}\(e^{V_{\n(t)}};~\frac{\mathbf{Z}_{\n(t)}}{t}\in \R^d\),
\end{equation}
то для того, чтобы изучить предельное поведение последовательности (\ref{n.1}), достаточно получить
предельные теоремы для последовательности
\begin{equation}{\label{n.1.2}}
\frac{1}{t}\ln {\bf E}\(e^{V_{\n(t)}};~\frac{\mathbf{Z}_{\n(t)}}{t}\in B\),
\end{equation}
что и осуществлено в теореме~\ref{th_th3}.
Следствием теоремы~\ref{th_th3},  в частном случае, когда $v=0$ п.н., является ПБУ для ОПВ $\mathbf{Z}(t)$ (возможно, обрывающегося), который сформулирован в  следствии~\ref{cor_cor2}.


\section{ Основные обозначения, основной результат}

Обозначим
$$
A(\lambda, \mm):=\ln \mathbf{E}\left(e^{\lambda \tau+\scalar{\mm}{\zz}+v} ; \tau<\infty\right),
$$
здесь и далее $\scalar{\cdot}{\cdot}$ --- скалярное произведение.
Функция $A(\lambda, \mm)$ является   преобразованием Лапласа над <<неполной>> мерой, отличающейся от меры
$$
{\bf E}\(e^v;~\t\in \cdot,~\zz\in \cdot\)=
{\bf E}\(e^v;~\t\in \cdot,~\zz\in \cdot,~\t<\infty\)+
{\bf E}\(e^v;~\zz\in \cdot,~\t=\infty\),
$$  отсутствием
второго слагаемого. Однако эта <<неполнота>>
полностью согласуется с тем очевидным обстоятельством, что характеристика
(\ref{n.1}), которая изучается, не зависит от распределения
$${\bf P}\(\zz\in \cdot,~v\in \cdot|~\t=\infty\).
$$

Рассмотрим два множества
$$
  {\cal A}^{\le 0}:=\{(\l,\mm):~A(\l,\mm)\le 0\},~~~
{\cal A}_\g^{\le 0}:=\{(\l,\mm):\l < \g,~A(\l,\mm)\le 0\},
$$
где $\g\in \overline{\R}:=\R\cup\{\infty\}$ --- фиксировано.
В качестве $\g$ будут выбираться:
\begin{itemize}
\item либо константы $\l_{\pm}$,
$0\le \l_+\le \l_-\le \infty$, где
\begin{equation}{\label{n.2}}
  \l_+:=\sup\{\l\ge 0:~{\bf E}e^{\l\t}<\infty\}=-\limsup_{t\to \infty}
  \frac{1}{t}\ln {\bf P}(\t>t),
\end{equation}
\begin{equation}{\label{n.3}}
  \l_-:=-\liminf_{t\to \infty}
  \frac{1}{t}\ln {\bf P}(\t>t);
 \end{equation}

\item либо константы $\l^*_{\pm}$,
$0\le \l^*_+\le \l^*_-\le \infty$, определенные при дополнительном условии
$$
  {\bf E}\(e^v\)<\infty
$$
соотношениями
\begin{equation*}{\label{n..2}}
  \l^*_+:=\sup\{\l\ge 0:~{\bf E}e^{v+\l\t}<\infty\}=-\limsup_{t\to \infty}
  \frac{1}{t}\ln {\bf E}(e^v;~\t>t),
\end{equation*}
\begin{equation*}{\label{n..3}}
  \l^*_-:=-\liminf_{t\to \infty}
  \frac{1}{t}\ln {\bf E}(e^v;~\t>t).
 \end{equation*}
\end{itemize}
Заметим, что если выполнено условие обрываемости $p={\bf P}(\t=\infty)>0$,
то выполняется $\l_+=\l_-=0$.

Построим новые функции
\begin{equation}\label{lem2_1.2}
 A(\mm):=-\sup\{\l:~(\l,\mm)\in {\cal A}^{\le 0}\},~~~
A_\g(\mm):=\max\{ -\gamma, A(\mm)\},
\end{equation}
где, по определению, считаем
$$
\sup\{\l:~\l\in \emptyset\}=-\infty.
$$
Для функции $H=H(\mm):~\R^d\to (-\infty, \infty]$
определим (см., например, \cite{Zalinescu2002}) преобразование Лежандра $H^{\mathfrak{Le}}=H^{\mathfrak{Le}}(\aa)$,
положив
$$
  H^{\mathfrak{Le}}(\aa):=\sup_{\mm}\{\scalar{\mm}{\aa}-H(\mm)\},~~~~\aa\in \R^d.
$$
Будем называть функцию $H=H(\aa)$,
отображающую $\R^{d}$  в $[0,\infty]$, {\it компактной}, если для любого
$c\ge 0$ множество $\skk{\aa:~H(\aa)\le c}$ есть компакт в  $\R^{d}$.
Легко показать, что любая компактная функция $H(\aa)$
полунепрерывна снизу.

Определим две функции: для $\aa\in \R^d$
$$
D(\aa):=A^{\mathfrak{Le}}(\aa),~~~ D_\g(\aa):=A_\g^{\mathfrak{Le}}(\aa).
$$

Следующая лемма будет доказана в разделе \ref{sec_proof_lem_1_2}.

\begin{lemma}\label{lem_th1}
\begin{itemize}
\item[\namedlabel{th1.1i}{$(i)$}] Функции $A(\mm)$, $A_{\gamma}(\mm)$ выпуклы и
полунепрерывны снизу.

\item[\namedlabel{th1.1ii}{$(ii)$}] Функции  $D(\aa)$, $D_{\gamma}(\aa)$ выпуклы,
полунепрерывны снизу и компактны.

\item[\namedlabel{th1.1iii}{$(iii)$}] Справедливы следующие формулы
\begin{equation}{\label{lem2_1.4}}
 A(\mm)=D^{\mathfrak{Le}}(\mm),~~~A_{\gamma}(\mm)=D_{\gamma}^{\mathfrak{Le}}(\mm),
\end{equation}
так что пары $A(\mm),D(\aa)$;  $A_{\gamma}(\mm),D_{\gamma}(\aa)$ являются парами взаимно
сопряженных (относительно преобразования Лежандра) функций.

\item[\namedlabel{th1.1iv}{$(iv)$}] Функции  $A_{\gamma}(\mm)$, $A(\mm)$ совпадают (и, следовательно, $D_{\gamma}(\aa)=D(\aa)$), тогда и только тогда,
когда выполнено условие
\begin{equation}{\label{lem2_1.5}}
\gamma \ge D(\0).
\end{equation}
\item[\namedlabel{th1.4ii}{$(v)$}]  Для всех $\aa\in \R^d$ справедливо
\begin{equation}{\label{1.11}}
D_{\gamma}(\aa)=\inf_{\h\in [0,1]}\skk{D(\h,\aa)+\gamma(1-\h)},
\end{equation}

где $$
D(\theta,\aa):= \sup\limits_{(\lambda,\mm) \in \mathcal{A}^{\le 0}}
 \{\lambda\theta + \scalar{\mm}{\aa}\}.
$$
\end{itemize}

\end{lemma}


Определим теперь две функции:
$$
D_+(\aa):=D_{\l_+}(\aa)-A_{\l_-}(0),~~~D_-(\aa):=D_{\l_-}(\aa)-A_{\l_+}(0),~~~\aa\in \R^d,
$$
где константы $\l_+$, $\l_-$ определены формулами (\ref{n.2}),
(\ref{n.3}), соответственно.
Заметим, что
$$
  A_{\l_\pm}(0)=\sup\{\l<\l_\pm:~{\bf E}e^{\l\t+v}\le 1\}.
$$
Заметим, что если выполнено условие обрываемости $p={\bf P}(\t=\infty)>0$,
то
$$
  D_+(\aa)=D_-(\aa)=D_0(\aa)-A_0(0),~~~~\aa\in \R^d.
$$
Для множества $B\in{\cal B}_d$ через $[B]$ и $(B)$
будем обозначать его замыкание и внутренность, соответственно.
Положим для $B\in{\cal B}_d$
$$
D_{\gamma}(B)=\inf\limits_{\aa\in B}D_{\gamma}(\aa).
$$

Основными результатами данной работы являются следующие утверждения.

\begin{theorem}\label{th_th3} Для любого борелевского  множества $B\subset \R^d$

                   \begin{equation}{\label{n.4}}
\limsup_{t\to \infty}\frac{1}{t}\ln {\bf E}
\(e^{V_{\n(t)}};~\frac{\mathbf{Z}(t)}{t}\in B\)\le -D_{\l_+}([B]),
\end{equation}
\begin{equation}{\label{n.5}}
 \liminf_{t\to \infty}\frac{1}{t}\ln {\bf E}
\(e^{V_{\n(t)}};~\frac{\mathbf{Z}(t)}{t}\in B\)\ge -D_{\l_-}((B)).
 \end{equation}
\end{theorem}
Используя равенство \eqref{16.12.1}, находим асимптотику вероятности \eqref{n.1}.

\begin{corollary}\label{cor_cor1} Для любого борелевского  множества $B\subset \R^d$

\begin{equation*}{\label{n.6}}
\limsup_{t\to \infty}\frac{1}{t}\ln {\bf P}_t
\(\frac{\mathbf{Z}((t)}{t}\in B\)\le -D_+([B]),
\end{equation*}
\begin{equation*}{\label{n.7}}
 \liminf_{t\to \infty}\frac{1}{t}\ln {\bf P}_t
\(\frac{\mathbf{Z}((t)}{t}\in B\)\ge -D_-((B)).
 \end{equation*}
\end{corollary}
Асимптотика вероятности \eqref{a04.11.1} очевидным образом извлекается из Теоремы~\ref{th_th3}, накладывая условие ${\bf P}(v=0)=1$.

\begin{corollary}\label{cor_cor2} Пусть ${\bf P}(v=0)=1$. Тогда для любого борелевского  множества $B\subset \R^d$

                   \begin{equation}{\label{n..4}}
\limsup_{t\to \infty}\frac{1}{t}\ln {\bf P}
\(\frac{\mathbf{Z}(t)}{t}\in B\)\le -D_{\l_+}([B]),
\end{equation}
\begin{equation}{\label{n..5}}
 \liminf_{t\to \infty}\frac{1}{t}\ln {\bf P}
\(\frac{\mathbf{Z}(t)}{t}\in B\)\ge -D_{\l_-}((B)).
 \end{equation}
В частности, если выполнено условие обрываемости $p={\bf P}(\t=\infty)>0$,
то $\l_-=\l_+=0$, и тогда выполнение неравенств (\ref{n..4}), (\ref{n..5})
означает, что семейство $\frac{\mathbf{Z}(t)}{t}$ удовлетворяет принципу
больших улонений в $\R^d$ с функцией уклонений $D_0(\aa)$.
\end{corollary}




\section{Доказательство Теоремы~\ref{th_th3}  }\label{sec_proof_th_3}
В основе доказательства Теоремы~\ref{th_th3} лежат следующие леммы, которые доказаны в разделе~\ref{sec_lem_proof}.

\begin{lemma}\label{lem_lem1} Для любого $\aa\in \R^d$
выполняется
\begin{equation}{\label{1.1}}
  \lim_{\v\downarrow 0}\limsup_{t\to \infty} \frac{1}{t}\ln
  {\bf E}\(e^{V_{\n(t)}};~\frac{{\mathbf{Z}}(t)}{t}\in (\aa)_\v\)\le -D_{\l_+}(\aa).
\end{equation}
\end{lemma}

\begin{lemma}\label{lem_lem2} Для любых $\aa\in \R^d$, $\v>0$
выполняется
\begin{equation}{\label{1.5}}
  \liminf_{t\to \infty} \frac{1}{t}\ln
  {\bf E}\(e^{V_{\n(t)}};~\frac{{\mathbf{Z}}(t)}{t}\in (\aa)_\v\)\ge -D(\aa),
\end{equation}
 \begin{equation}{\label{1.6}}
  \liminf_{t\to \infty} \frac{1}{t}\ln
  {\bf E}\(e^{V_{\n(t)}};~\frac{{\mathbf{Z}}(t)}{t}\in (\aa)_\v\)\ge -D_{\l_-}(\aa).
\end{equation}
\end{lemma}

Поскольку $A_{\l_-}(\mm)\ge A(\mm)$, то  $-D_{\l_-}(\aa)\ge -D(\aa)$ и,
следовательно,  неравенство (\ref{1.5}) следует из
неравенства (\ref{1.6}).   Однако
неравенство (\ref{1.5}), которое, вообще говоря, более грубое, чем
неравенство (\ref{1.6}), тем не менее представляет определенный интерес,
поскольку дает содержательную оценку снизу
в тех случаях, когда отсутствует информация о константе $\l_-$.

\begin{lemma}\label{lem_lem4} Для любого $N\in (0,\infty)$ найдется  $M\in (0,\infty)$
такое, что
$$
  \limsup_{t\to \infty}\frac{1}{t}\ln
  {\bf E}\(e^{V_{\n(t)}};~\frac{|\mathbf{Z}(t)|}{t}\ge M\)\le -N.
$$
\end{lemma}


Выполним теперь на базе лемм \ref{lem_lem1}--\ref{lem_lem4} доказательство Теоремы~\ref{th_th3}, которое повторяет основные шаги доказательства теоремы 4.1.1
\cite[стр. 259]{Bor1}.

\begin{proof}[Доказательство теоремы~\ref{th_th3}]

$(i)$. Оценка сверху (\ref{n.4}).
Фиксируем константы $\d>0$, $N<\infty$ и обозначим $D(B,\d,N):=\min\{N,D_{\l_+}([B])\}+\d,$
$$
L_+(B):=\limsup_{t\to \infty}\frac{1}{t}\ln
{\bf E}(e^{V_{\n(t)}};~\frac{\mathbf{Z}(t)}{t}\in B).
$$
В силу леммы \ref{lem_lem4} найдется компакт $K\subset \R^d$ такой, что для
$\overline{K}:=\R^d\setminus K$ выполняется
\begin{equation}{\label{1.26.11}}
L_+(\overline{K})\le -2D(B,\d,N).
\end{equation}
Далее, в силу леммы \ref{lem_lem1} для любого $\aa$ из компакта
$K\cap[B]$ найдется $\v(\aa)>0$ такое, что выполняется
\begin{equation}{\label{2.26.11}}
L_+((\aa)_{\v(\aa)})\le -D(B,\d,N).
\end{equation}
Получили открытое покрытие компакта $K\cap[B]$, из которого выделяем конечное подпокрытие
\begin{equation}{\label{3.26.11}}
  \{(\aa_i)_{\v(\aa_i)}\}_{i=1}^I:~~~~K\cap[B]\subset\cup_{i=1}^I(\aa_i)_{\v(\aa_i)}, ~~~
  I<\infty.
\end{equation}
Поэтому
$$
  L_+([B])\le L_+((K\cap[B])\cup \overline{K}),
  $$
и в силу (\ref{1.26.11}) -- (\ref{3.26.11})   имеем
$$
 L_+([B])\le -D(B,\d,N).
$$
Левая часть последнего неравенства   не зависит от $\d>0$
и $N<\infty$, поэтому неравенство сохранится, если в  его правой части
устремить  $\d\to 0$  и $N\to \infty$. Получили
оценка сверху (\ref{n.4}).

$(ii)$. Оценка снизу (\ref{n.5}).
Фиксируем $\aa\in (B)$, $\v>0$ таким образом, что $(\aa)_\v\subset (B)$ и обозначим
$$
L_-(B):=\liminf_{t\to \infty}\frac{1}{t}\ln
{\bf E}\left(e^{V_{\n(t)}};~\frac{\mathbf{Z}(t)}{t}\in B\right).
$$
Имеем
$$
L_-(B)\ge L_-((B))\ge L_-((\aa)_\v),
$$
поэтому, применяя лемму \ref{lem_lem2}, получаем
$$
L_-(B)\ge -D_{\l_-}(\aa).
$$
Левая часть последнего неравенства   не зависит от $\aa\in (B)$,
 поэтому неравенство сохранится, если   его правую часть
 максимизировать по $\aa\in (B)$. Получили
оценка снизу (\ref{n.5}).  Теорема 1 доказана.
\end{proof}


\section{Доказательство лемм~\ref{lem_lem1}, \ref{lem_lem2} и \ref{lem_lem4}}\label{sec_lem_proof}

Нам понадобится следующие обозначения. Для $(\theta,\aa) \in \mathbb{R}^{d+1}$ обозначим
\begin{equation}\label{1.6.001}
D_\Lambda(\h,\aa):=\inf_{r>0}r\Lambda\(\frac{\h}{r},\frac{\aa}{r}\),
\end{equation}
где
$$
  \Lambda(\h,\aa):=\sup_{{\l,\mm}}\skk{\l\h+\scalar{\mm}{\aa}-A(\l,\mm)},
  ~~~(\h,\aa)\in \R^{d+1}
$$
---{\it функция
уклонений}, которая определяется как
преобразование Лежандра функции $A(\l,\mm)$
\begin{equation*}
  \Lambda(\theta,\aa) = A^{\mathfrak{Le}}(\theta,\aa).
\end{equation*}
 Легко убедиться, что  функция $D_\Lambda(\h,\aa)$ выпукла и линейчата (т.е. линейна вдоль любого луча, выходящего из начала координат). Однако свойство полунепрерывности снизу для этой функции может отсутствовать.

В следующем утверждении мы приведем
дополнительные свойства функций $D(\aa)$, $D_{\gamma}(\aa)$:

\begin{lemma}\label{lem_th2}
\begin{itemize}
\item[\namedlabel{th1.4i}{$(i)$}]  Для всех $\aa\in \R^d$ справедливо

  \begin{equation}{\label{1.9}}
D(\aa)=\sup_{(\l,\mm)\in \mathcal{A}^{\le 0}}\skk{\l+\scalar{\mm}{\aa}};
\end{equation}

 \begin{equation}{\label{1.10}}
D_{\gamma}(\aa)=\sup_{(\l,\mm)\in \mathcal{A}^{\le 0}_{\gamma}}\skk{\l+\scalar{\mm}{\aa}};\ \
\end{equation}

\item[\namedlabel{th1.4iii}{$(ii)$}]   Для функции
$D_\Lambda(\h,\aa)$ (см. (\ref{1.6.001})
имеет место равенство

 \begin{equation}{\label{1.13}}
\( D_\Lambda^{\mathfrak{Le}}\)^{\mathfrak{Le}}(\h,\aa)=D(\h,\aa),
\end{equation}


и для всех $\theta > 0, \aa \in \mathbb{R}^d$ выполнено
 \begin{equation}{\label{1.14.2}}
D(\theta,\aa)=\lim_{\v \downarrow 0} \inf_{\aa'\in (\aa)_\v}D_{\Lambda}(\theta,\aa').
\end{equation}

\item[\namedlabel{th1.4iv}{$(iii)$}]  Если $\gamma<\infty$, то для функции
$$
\widehat{D}_{\gamma}(\aa):=\inf_{\h\in (0,1)}\skk{D(\h,\aa)+\gamma(1-\h)}
$$
имеет место равенство

 \begin{equation}{\label{1.15}}
\lim _{\varepsilon \downarrow 0} \inf _{\boldsymbol{\alpha}^{\prime} \in(\boldsymbol{\alpha}_\v)_{\boldsymbol{}}} \widehat{D}_{\gamma}\left(\boldsymbol{\alpha}^{\prime}\right)=D_{\gamma}(\aa).
\end{equation}
\end{itemize}
\end{lemma}

\begin{proof}[Доказательство леммы~\ref{lem_lem1}]
Из определения процесс $\mathbf{Z}(t)$ вытекает, что на событие $\skk{T_n\le t<T_n+\t_{n+1}}$ выполнено $\mathbf{Z}(t) = {\mathbf{Z}}_n.$ Следовательно,
\begin{equation}\label{lem1_1}
E(t):={\bf E}\(e^{V_{\n(t)}};~\frac{{\mathbf{Z}}(t)}{t}\in (\aa)_\v\)=\sum_{n\ge 0}
E_n(t),
\end{equation}
где
$$
  E_n(t):={\bf E}\(e^{V_{n}};~\frac{{\mathbf{Z}}_n}{t}\in (\aa)_\v,~T_n\le t<T_n+\t_{n+1}\).
$$

Оценим сначала часть суммы в \eqref{lem1_1} по $n>t^{2}$ :
$$
\sum_{n>t^{2}} E_n(t) \leqslant \sum_{n>t^{2}} {\bf E}\(e^{V_{n}};~T_n\le t\) =: \sum_{n>t^{2}} P_n(t).
$$
Выберем число $\l^*>0$ таким образом, что $A(-\l^*,\mathbf{0})\le 0$
(в силу условия
\ref{C} такая константа $\l^*$ всегда найдется), и рассмотрим
новые случайные  независимые  величины  $\t^*_j$   с распределением
$$
  {\bf P}(\t^*\in \cdot):=e^{-A(-\l^*,\mathbf{0})}{\bf E}\(e^{v-\l^*\t};~\t\in \cdot\).
$$
Легко показать, что $\p{\tau^* \in (0,\infty]} = 1,$ следовательно,  функция уклонений
$$
  \Lambda^*(\h):=\sup_{\l}\{\l\h-\ln {\bf E}e^{\l\t^*}\}
$$
неограниченно возрастает при монотонном приближении справа аргумента $\h$ к началу координат, т.е.
$ \lim\limits_{\h\downarrow 0}\Lambda^*(\h)=\infty. $
Далее, для $n\ge 1$ обозначим   $T^*_n:=\t^*_1+\cdots+\t^*_n,$
так что справедливо неравенство
\begin{equation*}
\begin{split}
P_n(t) & = e^{{\pm}n A(-\l^*,\mathbf{0})}{\bf E}\(e^{V_{n}\pm\l^*T_n };~T_n\le t\)  = e^{nA(-\l^*,\mathbf{0})}{\bf E}\(e^{\l^*T^*_n };~T^*_n\le t\) \\
&\le e^{nA(-\l^*,\mathbf{0})+\l^*t}{\bf E}\(T^*_n\le t\) \le e^{\l^*t}
  {\bf P}(T^*_n\le t),
\end{split}
\end{equation*}
где последнее неравенство справедливо в силу того, что $A(-\l^*,\mathbf{0})\leq0$.
Поэтому, в силу экспоненциального неравенства Чебышева при
$\frac{t}{n}\le {\bf E}\t^*$
имеем
$$
P_n(t)\le e^{\l^*t-n\Lambda^*\(\frac{t}{n}\)}.
$$
Таким образом, для $n\ge t^2$, $\frac{1}{t}\le {\bf E}\t^*$ имеем оценку
$$ P_n(t)\le e^{\l^*t}q^n(t), $$ где
$q(t):= e^{-\Lambda^*\(\frac{1}{t}\)} < 1 \text{ для всех достаточно больших } t.$ Следовательно,
$$
  \sum_{n\ge t^2}E_n(t)\le e^{\l^*t}\frac{q^{t^2}}{1-q}.
$$
  Поэтому
  \begin{equation}{\label{1.2}}
  \limsup_{t\to \infty} \frac{1}{t}\ln
  \sum_{n\ge t^2}E_n(t)=-\infty.
\end{equation}

Приведем теперь более точную оценку сверху для $E_n(t)$,
справедливую для всех $n\ge 1$. Для любого вектора
$(\l,\mm)\in \left({\cal A}^{\le 0}_{\l_+}\right)$      имеем
$$
  E_n(t)={\bf E}\(e^{V_n \pm\l{T}_n\pm\scalar{\mm}{\mathbf{Z}_n}};~
  \frac{{\mathbf{Z}}_n}{t}\in (\aa)_\v,~T_n\le t<T_n+\t_{n+1} \).
$$
На событии
$$
  \skk{\frac{{\mathbf{Z}}_n}{t}\in (\aa)_\v,~T_n\le t<T_n+\t_{n+1}}
$$
выполняется неравенство
$$
  e^{-\l{T_n}-\scalar{\mm}{\mathbf{Z}_n}} \le e^{-t(\l+\scalar{\mm}{\aa})+\sqrt{d}|\mm|\v t}
  \max\{1,~e^{\l\t_{n+1}} \}.
$$
Так как вектор $(\l,\mm)\in \left({\cal A}^{\le 0}_{\l_+}\right),$ то $\lambda < \lambda_+,$ и, следовательно
(учитывая, что $\lambda_+=0$ при $p=\mathbf{P}(\tau=\infty)>0$) будем иметь
$$
{\bf E} \max\skk{1,~e^{\l\t_{n+1}}}  < \infty.
$$
Таким образом, для $n\ge 1$ в любом случае имеем оценку
\begin{equation*}
\begin{split}
 E_n(t) & \le
 e^{-t(\l+\scalar{\mm}{\aa})+\sqrt{d}|\mm|\v t}
 {\bf E} \( \max\{1,~e^{\l\t_{n+1}}\} e^{V_n +\l{T}_n+\scalar{\mm}{\mathbf{Z}_n}} \) \\
 & \le  e^{-t(\l+\scalar{\mm}{\aa})+\sqrt{d}|\mm|\v t} {\bf E}
\max\skk{1,~e^{\l\t}}  e^{n A(\lambda,\mm)} \\
&  \le e^{-t(\l+\scalar{\mm}{\aa})+\sqrt{d}|\mm|\v t} {\bf E}
\max\skk{1,~e^{\l\t}}.
\end{split}
\end{equation*}
Из которой вытекает неравенство
\begin{equation}\label{lem1_2_2}
 \sum_{n = 1}^{[t^2]}E_n(t)\le  {\bf E}
\max\skk{1,~e^{\l\t}} t^2e^{-t(\l+\scalar{\mm}{\aa})+\sqrt{d}|\mm|\v t},
\end{equation}

Оценим наконец $E_0(t)={\bf P}\(\mathbf{0}\in (\aa)_\v,~\t>t\)$.
 Имеем
\begin{equation}\label{lem1_2_3}
 \lim_{\v\downarrow 0}\limsup_{t\to \infty} \frac{1}{t}\ln
E_0(t)=
\begin{cases} -\infty,~~~\mbox{если}~~~\aa\not = \mathbf{0};\\
          -\l_+,~~~\mbox{если}~~~\aa=\mathbf{0}.
\end{cases}
\end{equation}

Таким образом, из \eqref{lem1_2_2} и \eqref{lem1_2_3} получается неравенство
\begin{equation}\label{lem1_2}
 \limsup_{t\to \infty}\frac{1}{t}\ln
 \sum_{n = 0}^{[t^2]}E_n(t)\le -(\l+\scalar{\mm}{\aa})+\sqrt{d}|\mm|\v .
\end{equation}
Из соотношений \eqref{lem1_1}, \eqref{1.2}, \eqref{lem1_2} вытекает, что для любого $(\l,\mm)\in \left({\cal A}^{\le 0}_{\l_+}\right)$
\begin{equation}\label{lem1_3}
 \lim_{\v\downarrow 0}
 \limsup_{t\to \infty}\frac{1}{t}\ln E(t)\le -(\l+\scalar{\mm}{\aa}).
\end{equation}
Поскольку к любой точке $\left(\lambda, \mm\right)$ границы $\partial \mathcal{A}_{\l_+}^{\le 0}$ можно приблизиться точками $\left(\lambda_{n}, \mm_{n}\right)$ из внутренности $\left(\mathcal{A}_{\l_+}^{ \le 0}\right)$, то неравенство \eqref{lem1_3} справедливо для всех $(\lambda, \mm)$ из замкнутого выпуклого множества $\mathcal{A}_{\l_+}^{\le 0}$.
Минимизируя правую часть  неравенства \eqref{lem1_3} по $(\l,\mm)\in {\cal A}^{\le 0}_{\l_+}$
и используя утверждение \ref{th1.4i} леммы~\ref{lem_th2}, получаем
  \begin{equation*}{\label{1.3}}
 \lim_{\v\downarrow 0}
 \limsup_{t\to \infty}\frac{1}{t}\ln E(t) \le-D_{\l_+}(\aa).
\end{equation*}
Лемма~\ref{lem_lem1} доказана.
\end{proof}

\begin{proof}[Доказательство леммы~\ref{lem_lem2}]


Приведем сначала доказательство неравенства (\ref{1.5}). Выберем $\l^*>0$ такое, что выполняется $A(-\l^*,\mathbf{0})\le 0$,  что эквивалентно
{${\bf E}(e^{-\l^*\t+v})~\le~1$} (как уже отмечалось при
доказательстве леммы~\ref{lem_lem1}, в силу условия
\ref{C} такая константа $\l^*$ всегда найдется).
Для этого $\l^*$ построим случайный вектор $(\t^*,\zz^*,~v^*)\in \R_+\times\R^d\times\R$
с распределением
$$
    {\bf P}\(\t^*\in \cdot,~\zz^*\in \cdot,~v^*\in \cdot\):=
  $$
$$
   {\bf P}^*\(\t\in \cdot,~\zz\in \cdot,~v\in \cdot\):=
   {\bf E}\(e^{-\l^*\t+v};~\t\in \cdot,~\zz\in \cdot,~v\in \cdot\).
$$
Заметим, что в случае, когда $A(-\l^*,\mathbf{0})<0$,
распределение этого вектора будет несобственным, т.е.
$$
{\bf P}\(\t^*\in (0,\infty),~\zz^*\in \R^d,~v^*\in \R\)=e^{A(-\l^*,\mathbf{0})}<1.
$$
Это распределение можно произвольным образом доопределить на множестве $\skk{\infty}\times\mathbb{R}^d\times\mathbb{R}.$
Преобразование Лапласа над распределением вектора $(\t^*,\zz^*)$ обозначим как
$$
e^{A^*(\l,\mm)}:={\bf E}e^{\l\t^*+\scalar{\mm}{\zz^*}}=
{\bf E}\(e^{\l\t+\scalar{\mm}{\zz}-\l^*\t+v}\)=e^{A(-\l^*+\l,\mm)},
$$
т.е. положим
$$
A^*(\l,\mm):=A(-\l^*+\l,\mm).
$$

Определим далее последовательность
$\{(\t^*_i,\zz^*_i,v^*_i);~i=1,\dots\}$
независимых копий случайного вектора $(\t^*,\zz^*,v^*)$ и для $n=0,1,\dots$ обозначим через
$\(T^*_n,\mathbf{Z}^*_n,V^*_n\)$
 частичные суммы этих векторов. Процесс восстановления определим естественным образом
$\n^*(t):=\sup\{n\ge 0:~T^*_n\le t\}.$
Тогда можно определить новую пару ОПВ:
$$
  \({T}^*(t),{\mathbf{Z}}^*(t)\):=\({T}^*_{\n^*(t)},{\mathbf{Z}}^*_{\n^*(t)}\).
$$
Поясним, как определяется распределение (вообще говоря, несобственное)
этой новой пары  ОПВ:
\begin{equation*}
\begin{split}
 & {\bf P}\({\mathbf{Z}}^*(t)\in \cdot,~T^*(t)\in \cdot\)  =
  \sum_{n\ge 0}{\bf P}\({\mathbf{Z}}^*_n\in \cdot,~T^*_n\in \cdot;~
  T^*_n\le t<T^*_n+\t^*_{n+1}\) \\
& =
\sum_{n\ge 0}{\bf E}\(e^{-\l^*T_{n+1}+V_{n+1}};~{\mathbf{Z}}_n\in \cdot;~
T_n\in \cdot,~T_n\le t<T_n+\t_{n+1}\).
\end{split}
\end{equation*}

%

Для $i \geq 1$ обозначим $\overline{v}_i:=-\l^*\t_i+v_i,$
так что выполняется
$$
  {\bf P}^*\(\t_i\in \cdot,~\zz_i\in \cdot,~v_i\in \cdot\)
  :=
  {\bf E}\(e^{\overline{v}_i};~\t_i\in \cdot,~\zz_i\in \cdot,~v_i\in \cdot\).
$$
Не трудно видеть, что найдутся такие константы
$$
q>0,~~0<c<\infty,~~0<R<\infty,
$$
что  для событий $ B_i:=\{c<\t_i,~~|\zz_i|\le R,~~|v_i|\le R\}, i \ge 1$
выполняются неравенства
\begin{equation*}{\label{1.7..}}
  {\bf P}^*(B_i)\ge q.
\end{equation*}
Для $T>0$ обозначим  $k(T):=\frac{T}{c}$, 
$$
 B(T):=\{c<\t_i,~~|\zz_i|\le R,~~
 |v_i|\le R,~~\mbox{для всех}~~i\le \n(T)+1\}.
 $$
\begin{lemma}\label{lem_lem3} Для любого $T>0$
$$
{\bf P}^*\(B(T)\)\ge q^{k(T)+1}.
$$
\end{lemma}

\begin{proof}
Достаточно заметить, что если $\t_i>c$, $i\ge 1$, то
 справедливо включение
\begin{equation*}
    \bigcap\limits_{i=1}^{[k(T)]+1} B_i \subseteq B(T).
\end{equation*}
Лемма \ref{lem_lem3} доказана.
\end{proof}

Продолжим доказательство (\ref{1.5}).  Очевидно, что
\begin{equation*}
\begin{split}
 &  {\bf E}\(e^{V_{\n(t)}};~\frac{{\mathbf{Z}}_{\n(t)}}{t}\in (\aa)_{2\v}\)=
  {\bf E}\(e^{\l^* T_{\n(t)}+\overline{V}_{\n(t)}};
  ~\frac{{\mathbf{Z}}_{\n(t)}}{t}\in (\aa)_{2\v}\) \\
& \ge \int_{T=0}^{2\d t}
{\bf E}\(e^{\l^* T_n+\overline{V}_n};
 ~\frac{{\mathbf{Z}}_n}{t}\in (\aa)_{\v},~~\frac{{T}_n}{t}\in (1-\d)_{\d},~~t-T_n\in dT\)
  \times I(T,\v),
  \end{split}
\end{equation*}
где
$$
 I(T,\v):=
 {\bf E}\(e^{\l^* T_{\n(T)}-\overline{v}_{\n(T)+1}+\overline{V}_{\n(T)+1}};
  ~\frac{{\mathbf{Z}}_{\n(T)}}{t}\in (\mathbf{0})_{\v}\cap B(T)\).
$$
Поскольку на событии $\skk{\frac{T_n}{n}\in \frac{t}{n}(1-\d)_\d}$
при $T\le 2\d t$ выполняется
$$
\l^*T_n\ge \l^*t(1-2\d)
$$ и на событии $B(T)$
при $T\le 2\d t$ и $2R\d\le c\v$
выполняются соотношения
$$\l^*T_{\n(T)}-\overline{v}_{\n(T)+1}\ge -R,~~
\skk{\frac{{\mathbf{Z}}_{\n(T)}}{t}\in (\mathbf{0})_{\v}}\cap B(T)=B(T),
$$
то имеем
\begin{equation*}
    \begin{split}
& e^{-\l^*t(1-2\d)}{\bf E}\(e^{V_{\n(t)}};~\frac{{\mathbf{Z}}_{\n(t)}}{t}\in (\aa)_{2\v}\) \\
&  \ge \int_{T=0}^{2\d t}
  {\bf E}\(e^{\overline{V}_n};
  ~\frac{{\mathbf{Z}}_n}{n}\in \frac{t}{n}(\aa)_{\v},~
  \frac{T_n}{n}\in \frac{t}{n}(1-\d)_\d,~t-T_n\in dT\)\times
  e^{-R}J(T,\v),
\end{split}
\end{equation*}
  где
  $$
   J(T,\v):=
 {\bf E}\(e^{\overline{V}_{\n(T)+1}};
  ~B(T)\)={\bf P}^*(B(T))\ge q^{k(T)+1}.
  $$
В последнем неравенстве мы воспользовались леммой \ref{lem_lem3}.
Мы получили
\begin{equation}{\label{1.7}}
\begin{split}
& {\bf E}\(e^{V_{\n(t)}};~\frac{{\mathbf{Z}}_{\n(t)}}{t}\in (\aa)_{2\v}\) \\
& \ge e^{\l^*t(1-2\d)-R}q^{\frac{2\d t}{c}+1}\int_0^{2t\d}
{\bf P}^*\(
  ~\frac{{\mathbf{Z}}_n}{n}\in \frac{t}{n}(\aa)_{\v},~
  \frac{T_n}{n}\in \frac{t}{n}(1-\d)_\d,~t-T_n\in dT\) \\
& = e^{\l^*t(1-2\d)-R}q^{\frac{2\d t}{c}+1}
{\bf P}^*\(
  ~\frac{{\mathbf{Z}}_n}{n}\in \frac{t}{n}(\aa)_{\v},~
   \frac{T_n}{n}\in \frac{t}{n}(1-\d)_\d\).
   \end{split}
\end{equation}

Чтобы продолжить доказательство формулы (\ref{1.5}), нам понадобится следующее утверждение.

\begin{lemma}\label{lem_lem3*} Для любых $\v>0,~ \d>0,~r>0$ имеет место следующая
оценка снизу в принципе больших уклонений для сумм
$(\frac{T_n}{n},\frac{\mathbf{Z}_n}{n})$; $n:=[rt]$,
для несобственного, вообще говоря, распределения ${\bf P}^*(\cdot)$:
\begin{equation}{\label{1.1n}}
\liminf_{t\to \infty}\frac{1}{t}\ln
{\bf P}^*\left(\left(\frac{T_n}{n},\frac{\mathbf{Z}_n}{n}\right)\in \frac{t}{n}(1-\d)_\d\times(\aa)_\v\right)\ge
-\Lambda^*_r((1-\d)_\d\times(\aa)_\v),
\end{equation}
где $\Lambda^*_r(\h,\aa):=r\Lambda^*(\frac{(\h,\aa)}{r}) $ и где для множества
$B\subset \R_+\times \R^d$
$$
\Lambda^*_r(B):=\inf_{(\h,\aa)\in B}\Lambda^*_r(\h,\aa).
$$
\end{lemma}
\begin{proof}[Доказательство леммы~\ref{lem_lem3*}]
 Рассмотрим наряду с несобственным распределением
${\bf P}^*(\t\in \cdot,\zz\in \cdot)$
собственное распределение
\begin{equation}{\label{1.4n}}
 \hat{{\bf P}}(\t\in \cdot,\zz\in \cdot):=e^{-C}{\bf P}^*(\t\in \cdot,\zz\in \cdot),
\end{equation}
где  $C:=\ln{\bf E}e^{\overline{v}}$. Функцию уклонений, отвечающую
$\hat{{\bf P}}$--распределению вектора $(\t,\zz)$ обозначим
$$
  \hat{\Lambda}(\h,\aa):=\sup_{(\l,\mm)}\{\l\h+\scalar{\mm}{\aa}-\ln\hat{\bf E}e^{\l\t+\scalar{\mm}{\zz}} \}.
  $$
  Очевидно, что справедливо равенство
\begin{equation}{\label{1.3n}}
\hat{\Lambda}(\h,\aa)=\Lambda^*(\h,\aa)+C.
\end{equation}
Воспользуемся теперь известной (см., например, \cite[теорема 1.2.1]{Bor1})
оценкой снизу в принципе больших уклонений для сумм
$(\frac{T_n}{n},\frac{\mathbf{Z}_n}{n})$; $n:=[rt]$,
для собственного распределения $\hat{{\bf P}}(\cdot)$:
\begin{equation}{\label{1.2n}}
\liminf_{t\to \infty}\frac{1}{t}\ln
\hat{{\bf P}}\left(\left(\frac{T_n}{n},\frac{\mathbf{Z}_n}{n}\right)\in \frac{t}{n}(1-\d)_\d\times(\aa)_\v\right)
\ge
-\hat{\Lambda}_r((1-\d)_\d\times(\aa)_\v),
\end{equation}
где $\hat{\Lambda}_r(\h,\aa):=r\hat{\Lambda}\left(\frac{(\h,\aa)}{r}\right) $ и где для множества
$B\subset \R_+\times \R^d$
$$
\hat{\Lambda}_r(B):=\inf_{(\h,\aa)\in B}\hat{\Lambda}_r(\h,\aa).
$$
Остается заметить, что в силу (\ref{1.4n})  и (\ref{1.3n}) левая (правая) часть
(\ref{1.1n}) отличается от  левой (правой) части (\ref{1.2n}) на слагаемое $-rC$.
\end{proof}

Продолжим доказательство неравенства (\ref{1.5}).
Используя (\ref{1.7}) и лемму \ref{lem_lem3*}, получаем
\begin{equation}{\label{08.12.1.9}}
\begin{split}
  L_-(\aa, 2\v) &:=\liminf_{t\to \infty}
  \frac{1}{t}\ln {\bf  E}\left(e^{V_{\n(t)}}:~\frac{Z_{\n(t)}}{t}\in (\aa)_{2\v}\right) \\
  & \ge  -\Lambda^*_r((1-\d)_\d\times (\aa)_\v)+\l^*-W\d,
  \end{split}
\end{equation}
где $W:=(4\l^*+ \frac{4}{c}|\ln q|)$. Максимизируя правую часть (\ref{08.12.1.9}) по
$r>0$, используя обозначения
\begin{equation*}
\begin{split}
&  D^*_{\Lambda^*}(\h,\bb):=\inf_{r>0}\Lambda^*_r(\h,\bb),~~~\h>0,~\bb\in \R^d, \\
&  D^*_{\Lambda^*}(B):=\inf_{(\h,\bb)\in B}D^*_{\Lambda^*}(\h,\bb),~~~B\subset
  (0,\infty)\times \R^d,
 \end{split}
\end{equation*}
и равенство
$$
  \inf_{r>0}\Lambda^*_r((1-\d)_\d\times (\aa)_\v)=
  D^*_{\Lambda^*}((1-\d)_\d\times (\aa)_\v),
$$
получаем
$$
L_-(\aa, 2\v)\ge -D^*_{\Lambda^*}((1-\d)_\d\times (\aa)_\v)+\l^*-W\d.
$$
Поскольку для любого $u>0$ выполняется
$$
   D^*_{\Lambda^*}(u\h,u\bb)=uD^*_{\Lambda^*}(\h,\bb),
$$
то из последнего неравенства выводим
\begin{equation}{\label{1.9.}}
\begin{split}
L_-(\aa, 2\v) & \ge -(1-\d)D^*_{\Lambda^*}((1)_{\d'}\times (\bb)_{\v'})+\l^*-W\d
 \\
& \ge -(1-\d)D^*_{\Lambda^*}(\{1\}\times (\bb)_{\v'})+\l^*-W\d,
\end{split}
\end{equation}
где $\d':=\frac{\d}{1-\d}$, $\v':=\frac{\v}{1-\d}$, $\bb:=\frac{\aa}{1-\d}$.
Заметим, что для любого $\v>0$ найдется $\d_0=\d_0(\v)>0$ такое, что для всех
$\d\in (0,\d_0)$ выполняется
$$
 (\aa)_{\v/2}\subset (\bb)_{\v'},
$$
и, следовательно,
$$
-D^*_{\Lambda^*}(\{1\}\times (\bb)_{\v'})\ge -D^*_{\Lambda^*}(\{1\}\times (\aa)_{\v/2}).
$$
Используя  последнее неравенство  для оценки снизу  правой части (\ref{1.9.}),
получаем
$$
L_-(\aa, 2\v)\ge -(1-\d)D^*_{\Lambda^*}(\{1\}\times (\aa)_{\v/2})+\l^*-W\d;
$$
устремляя $\d\downarrow 0$, имеем для любого $N\ge 2$
$$
L_-(\aa, 2\v)\ge -D^*_{\Lambda^*}(\{1\}\times (\aa)_{\v/2})+\l^*\ge
-D^*_{\Lambda^*}(\{1\}\times (\aa)_{\v/N})+\l^*;
$$
устремляя $N\to \infty$  и используя (4.5), получаем неравенство
\begin{equation}{\label{1.17}}
L_-(\aa, 2\v)\ge
-(D^*_{\Lambda^*}(1,\aa)+\l^*).
\end{equation}

Нам осталось установить взаимосвязь между функциями $D^*(u,\aa)$ и $D(\aa).$  Для любого $0<u\le 1$  воспользуемся представлениями
$$
  D^*(u,\aa)=\sup_{A^*(\l,\mm)\le 0}\{{\l}u+\scalar{\mm}{\aa}\},~~~
  A^*(\l,\mm)=A(\l-\l^*,\mm),
$$
в силу которых получаем равенство
\begin{equation}{\label{1.17.u}}
 D^*(u,\aa)-{\l^*}u=\sup_{A(\l-\l^*,\mm)\le 0}\{\l{u}+\scalar{\mm}{\aa}\}-\l^*{u}=D(u,\aa).
\end{equation}

Применяя (\ref{1.17.u})  при $u=1$ к неравенству (\ref{1.17}) и используя (\ref{08.12.1.9}),
получаем доказательство неравенства (\ref{1.5}). При этом получено доказательство
неравенства (\ref{1.6}) в случае, когда выполнено
\begin{equation}{\label{4.26.11}}
D(\aa)=D_{\l_-}(\aa).
\end{equation}
Докажем теперь
неравенство (\ref{1.6}) в случае, когда условие (\ref{4.26.11})  не выполнено, т.е. когда
\begin{equation}{\label{5.26.11}}
  D(\aa)>D_{\l_-}(\aa).
\end{equation}
Заметим, что в случае  (\ref{5.26.11}) выполняется $\l_{-}<D(\mathbf{0})$
(см. лемму \ref{lem_th1}, пункт \ref{th1.1iv}) и, следовательно,
\begin{equation}{\label{7.26.11}}
\l_-<\infty.
\end{equation}
Поэтому нам достаточно доказать
неравенство (\ref{1.6}) в случае (\ref{7.26.11}). Выполним это доказательство.
В этом случае последний <<большой скачок>>  $\t_{\nu(t) + 1}$ вносит некоторый вклад в асимптотику исследуемой вероятности.
 Очевидно, что
\begin{equation*}
\begin{split}
  & {\bf E}\(e^{V_{\n(t)}};~\frac{{\mathbf{Z}}_{\n(t)}}{t}\in (\aa)_{2\v}\)=
  {\bf E}\(e^{\l^* T_{\n(t)}+\overline{V}_{\n(t)}};
  ~\frac{{\mathbf{Z}}_{\n(t)}}{t}\in (\aa)_{2\v}\)\\
   &\ge \  {\bf E}\(e^{\l^* T_n+\overline{V}_n};
 ~\frac{{\mathbf{Z}}_n}{t}\in (\aa)_{\v},~~\frac{{T}_n}{t}\in (u-\d)_{\d}\)
  \times {\bf P}(\t>t(1-u+2\d)),
  \end{split}
\end{equation*}
где  число $u\in (0,1)$ фиксировано. Поскольку на событии
$\{\frac{T_n}{n}\in \frac{t}{n}(u-\d)_\d\}$
 выполняется
$$
\l^*T_n\ge \l^*t(u-2\d),
$$ то имеем
\begin{equation}{\label{1.18}}
\begin{split}
& {\bf E}\(e^{V_{\n(t)}};~\frac{{\mathbf{Z}}_{\n(t)}}{t}\in (\aa)_{2\v}\) \\
& \ge \ e^{\l^*t(u-2\d)}  {\bf E}\(e^{\overline{V}_n};
  ~\frac{{\mathbf{Z}}_n}{n}\in \frac{t}{n}(\aa)_{\v},
  \frac{T_n}{n}\in \frac{t}{n}(u-\d)_\d\)\!\times\!
  {\bf P}(\t>t(1-u+2\d)).
  \end{split}
   \end{equation}
Далее, повторяя с очевидными изменениями вывод из
(\ref{1.7}) неравенства  (\ref{1.17}), выводим из (\ref{1.18}) для всех $u \in (0,1), \aa \in \mathbb{R}^d$
неравенство
\begin{equation}{\label{1.19}}
   L_-(\aa,2\v)\ge -(D^*(u,\aa)-\l^*{u}+\l_-(1-u)).
\end{equation}
Применяя  (\ref{1.17.u}) к правой части неравенства (\ref{1.19}),
получаем
$$
  L_-(\aa,2\v)\ge -(D(u,\aa)+\l_-(1-u)).
$$
Максимизируя правую часть последнего неравенства по $u\in (0,1)$,
получаем для всех $\aa \in \mathbb{R}^d$
\begin{equation}\label{lem32_add1}
  L_-(\aa,2\v)\ge -\widehat{D}_{\lambda_-}(\aa),
\end{equation}
где функция $\widehat{D}_{\lambda_-}(\aa)$ определена в \ref{th1.4iv} леммы~\ref{lem_th2}.
Выберем теперь произвольные $\aa^{\prime} \in(\aa)_{\varepsilon}$ и $\varepsilon^{\prime}>0$ такие, что выполняется $\left(\aa^{\prime}\right)_{\varepsilon^{\prime}} \subset(\aa)_{\varepsilon} .$ Применяя \eqref{lem32_add1} для $\aa^{\prime}$ и $\varepsilon^{\prime}$, получаем
\begin{equation}\label{lem32_add2}
     L_-(\aa,2\v)\ge  L_-(\aa^{\prime},2\v^{\prime}) \ge -\widehat{D}_{\lambda_-}(\aa^{\prime}),
\end{equation}
Максимизируя далее правую часть \eqref{lem32_add2} по  $\aa^{\prime} \in(\aa)_{\varepsilon}$, получаем для любого $\varepsilon^{\prime} \in (0, \varepsilon]$
$$
 L_-(\aa,2\v)\ge -\inf _{\aa^{\prime} \in(\aa)_{e^{\prime}}}  \widehat{D}_{\lambda_-}(\aa^{\prime}).
$$
Осталось воспользоваться равенством \eqref{1.15} и получить утверждение (\ref{1.6}) при дополнительном условии $\l_-<\infty$.
 Лемма~\ref{lem_lem2} доказана.
\end{proof}

\begin{proof}[Доказательство леммы~\ref{lem_lem4}] Легко видеть, что для любых $\gamma>0$, $\widetilde{\lambda}>0$ п.н. справедливы неравенства
\begin{equation} \label{28.06.1}
\mathbf{I}\left(\frac{|\mathbf{Z}(t)|}{t}\geq M\right)\leq\mathbf{I}\left(\gamma|\mathbf{Z}(t)|-\widetilde{\lambda} T_{\nu(t)}\geq M\gamma t-\widetilde{\lambda} t\right)
\leq\frac{e^{\gamma|\mathbf{Z}_{\nu(t)}|-\widetilde{\lambda} T_{\nu(t)}}}{e^{M\gamma t-\widetilde{\lambda} t}},
\end{equation}
\begin{equation} \label{28.06.2}
\mathbf{I}(\nu(t)=k)\leq\mathbf{I}(T_k\leq t)=\mathbf{I}(e^{-T_k}\geq e^{-t})\leq\frac{e^{-T_k}}{e^{-t}}.
\end{equation}
Из условия $[\mathbf{C}^*]$ следует, что найдутся $\gamma>0$ и $\widetilde{\lambda}>0$ такие, что
\begin{equation} \label{28.06.3}
u:=\mathbf{E}e^{v+\gamma|\zz|-(\widetilde{\lambda}+1) \tau}<1.
\end{equation}
Выбирая $\gamma>0$ и $\widetilde{\lambda}>0$ так, чтобы было выполнено неравенство
(\ref{28.06.3}), используя неравенства (\ref{28.06.1}), (\ref{28.06.2}) и лемму Беппо Леви, получаем
\begin{equation}\label{28.06.5}
\begin{split}
\mathbf{E}\left(e^{V_{\nu(t)}};\frac{|\mathbf{Z}(t)|}{t}  \geq M\right) & \leq\mathbf{E} \left(\frac{e^{V_{\nu(t)}+\gamma|\mathbf{Z}_{\nu(t)}|-\widetilde{\lambda} T_{\nu(t)}}}{e^{M\gamma t-\widetilde{\lambda} t}}\right) \\
& \leq e^{-M\gamma t+\widetilde{\lambda} t}+\sum\limits_{k=1}^{\infty}\mathbf{E}\left( \frac{e^{V_{k}+\gamma|\mathbf{Z}_{k}|-\widetilde{\lambda} T_{k}}}{e^{M\gamma t-\widetilde{\lambda} t}}; \ \nu(t)=k \right) \\
& \leq e^{-M\gamma t+\widetilde{\lambda} t}+ \sum\limits_{k=1}^{\infty}\mathbf{E}\left( \frac{e^{V_{k}+\gamma|\mathbf{Z}_{k}|-(\widetilde{\lambda}+1) T_{k}}}{e^{M\gamma t-(\widetilde{\lambda}+1) t}} \right)\\
& \leq e^{-M\gamma t+\widetilde{\lambda} t}+e^{-t(M\gamma-\widetilde{\lambda}-1)}
\sum\limits_{k=1}^{\infty}\left(\mathbf{E}e^{v+\gamma|\zz|-(\widetilde{\lambda}+1) \tau}\right)^k \\
& \leq\left(1+\frac{u}{1-u}\right)e^{-t(M\gamma-\widetilde{\lambda}-1)}.
\end{split}
\end{equation}
Используя неравенство (\ref{28.06.5}), выбирая $M=\frac{N+\widetilde{\lambda}+1}{\gamma}$, будем иметь
$$
\limsup\limits_{t\rightarrow\infty}\frac{1}{t}\ln \mathbf{E}\left(e^{V_{\nu(t)}};\frac{|\mathbf{Z}(t)|}{t}\geq M\right)
\leq-M\gamma+\widetilde{\lambda}+1= -N.
$$
\end{proof}

\section{Доказательство свойств основных функций (леммы~\ref{lem_th1} и \ref{lem_th2})}\label{sec_proof_lem_1_2}

Доказательство лемм ~\ref{lem_th1} и \ref{lem_th2} в значительной степени  повторяет доказательство аналогичных лемм в работе \cite{MogPro2019}. Однако, для удобства читателя мы приводим полные доказательства. 

\subsection{Свойства преобразования Лежандра}\label{subs2.1}
В дальнейшем нам понадобятся
некоторые свойства выпуклых полунепрерывных снизу
функций $F=F(\boldsymbol{u})$, отображающих $\R^d$  в $(-\infty,\infty]$, и преобразований Лежандра над ними.
Обозначим класс таких функций через $\mathcal{F}=\mathcal{F}_d$.
Известны следующие свойства преобразования Лежандра (см., например, \cite{BarbuPrecupanu2012} или \cite{Zalinescu2002}):
\begin{itemize}
\item[\namedlabel{L1}{$\mathbf{[L_1]}$}]  {\it Для любой функции $F\in \mathcal{F}$
справедливо $F^\mathfrak{Le}\in \mathcal{F}$,
т.е. выпуклая полунепрерывная снизу
функция переводится преобразованием Лежандра в выпуклую полунепрерывную
 снизу функцию} (см., например, Теорему 2.3.1 (i), стр. 75 в \cite{Zalinescu2002} или Предложение 2.19, (i), стр. 77 в \cite{BarbuPrecupanu2012}).

\item[\namedlabel{L2}{$[\mathbf{L_2}]$}] {\it Для любой функции $F\in \mathcal{F}$ справедливо
$$
   \(F^\mathfrak{Le}\)^\mathfrak{Le}=F,
$$
т.е. повторное применение преобразование Лежандра переводит выпуклую
полунепрерывную снизу функцию в себя} (см., например, Теорему 2.3.1 (iv), стр. 75 в \cite{Zalinescu2002} или Предложение 2.19, (iii), стр. 77 в \cite{BarbuPrecupanu2012}).
\end{itemize}
Для двух функций $F_1,~F_2\in \mathcal{F}$
определим операцию  свертки $*$, положив
$$
  F_1*F_2(\boldsymbol{u}):=\inf_{\boldsymbol{v}}\{F_1(\boldsymbol{v})+F_2(\boldsymbol{u}-\boldsymbol{v})\}.
$$
Нам понадобится следующее свойство преобразования Лежандра (см., например, Теорему 2.3.1 (ix), стр. 76 в \cite{Zalinescu2002})
\begin{itemize}
\item[\namedlabel{L3}{$[\mathbf{L_3}]$}] {\it Для любых $G_1,~G_2,~F_1,~F_2\in \mathcal{F}$
справедливо
}
$$
  (G_1*G_2)^\mathfrak{Le}=G_1^\mathfrak{Le}+G_2^\mathfrak{Le},~~~~  F^\mathfrak{Le}_1*F^\mathfrak{Le}_2=(F_1+F_2)^\mathfrak{Le} .
$$
\end{itemize}

Для произвольной выпуклой функции $F=F(\boldsymbol{u})$, отображающей
$\R^d$ в $(-\infty,\infty]$, через $cl\,F=cl\,F(\boldsymbol{u})$
обозначим
 наибольшую функцию из класса $\mathcal{F}$, минорирующую $F$.
Иначе говоря: (1) функция $cl\,F$ принадлежит классу $\mathcal{F}$;
(2) она во всех точках $\boldsymbol{u}$ не превышает $F(\boldsymbol{u})$; (3) для любой функции
$G(\boldsymbol{u})\in \mathcal{F}$ из неравенства
$$
  F(\boldsymbol{u})\ge G(\boldsymbol{u})~~~\mbox{для всех}~~~\boldsymbol{u}\in \R^d
$$
следует неравенство
$$
  cl\,F(\boldsymbol{u})\ge G(\boldsymbol{u})~~~\mbox{для всех}~~~\boldsymbol{u}\in \R^d.
$$
Следующее свойство хорошо известно (см., например, Теорему 2.3.4 (i), стр. 78 в \cite{Zalinescu2002})
\begin{itemize}
\item[\namedlabel{L4}{$[\mathbf{L_4}]$}] {\it Для любой выпуклой функции
$F=F(\boldsymbol{u})$: $\R^d\to (-\infty,\infty]$,
справедливо
$$
 \(F^\mathfrak{Le}\)^\mathfrak{Le}=cl\,F.
$$
}
Нам понадобится еще одно свойство

\item[\namedlabel{L5}{$[\mathbf{L_5}]$}] {\it Для любой выпуклой функции
$F=F(\boldsymbol{u})$: $\R^d\to (-\infty,\infty]$,
справедливо}
\begin{equation}{\label{2.1}}
   \lim_{\v\downarrow 0}\inf_{\boldsymbol{v}\in (\boldsymbol{u})_\v}F(\boldsymbol{v})=
   \(F^\mathfrak{Le}\)^\mathfrak{Le}(\boldsymbol{u}).
  \end{equation}
Свойство \ref{L5} используется в литературе, как факт не требующий доказательства (см., например, формулу (2.9) и рассуждения перед ней, стр. 71 в в \cite{BarbuPrecupanu2012}). Однако, доказательство этого свойства можно найти в работе \cite{MogPro2019}
\end{itemize}

\subsection{ Доказательство леммы~\ref{lem_th1}}\label{subs2.2}
\ref{th1.1i}.  Функция $A(\l,\mm)$ при любом
фиксированном $\mm\in \R^d$ неограниченно возрастает с ростом $\l$. Поэтому,
учитывая определение (\ref{lem2_1.2}), получаем свойство:
 \begin{equation}{\label{2.5}}
A(\mm)>-\infty~~~\mbox{{\it для любого}}~~~\mm\in \R^d.
\end{equation}
Пусть, далее, $A(\mm)<\infty$, и последовательность $\l_n$ такова, что
$$
  A(\l_n,\mm)\le 0~~~\mbox{при всех}~~~n\ge 1~~~\mbox{и}~~~
 \lim_{n\to \infty}\l_n=-A(\m).
$$
Тогда, в силу полунепрерывности снизу функции $A(\l,\mm)$,
получаем $$A(-A(\mm),\mm)~\le~0.$$ Другими словами, мы установили свойство:
 \begin{equation}{\label{2.6}}
\mbox{{\it из}}~~~A(\mm)<\infty~~~\mbox{{\it следует}}~~~
A(-A(\mm),\mm)\le 0.
 \end{equation}

Убедимся теперь, что функция $A(\mm)$ выпукла: {\it для любых
$\mm_1,\mm_2\in \R$, $p, q~\ge~0$, $p+q=1$ выполняется}
\begin{equation}{\label{2.7}}
A(p\mm_1+q\mm_2)\le pA(\mm_1)+qA(\mm_2).
 \end{equation}
Если $A(\mm_1)=\infty$, или $A(\mm_2)=\infty$, то неравенство
(\ref{2.7}) выполнено. Если $A(\mm_1)<\infty$ и $A(\mm_2)<\infty$,
то в силу (\ref{2.6})     имеем
$$
 A(-A(\mm_1),\mm_1)\le 0,~~~A(-A(\mm_2),\mm_2)\le 0.
$$
Поэтому,
в силу выпуклости функции $A(\l,\mm)$ справедливо
$$
 A(-pA(\mm_1)-qA(\mm_2),p\mm_1+q\mm_2)\le
 pA(-A(\mm_1),\mm_1)+qA(-A(\mm_2),\mm_2)\le 0,
$$
т.е.
$$
   A(-pA(\mm_1)-qA(\mm_2),p\mm_1+q\mm_2)\le 0.
$$
Из последнего следует
$-A(p\mm_1+q\mm_2)\ge -pA(\mm_1)-qA(\mm_2)$, т.е.
 (\ref{2.7}). Выпуклость функции $A(\mm)$
установлена.

Убедимся теперь, что функция $A(\mm)$ полунепрерывна снизу:
для любой последовательности $\mm_n$, сходящейся к $\mm$ при
$n\to \infty$
\begin{equation}{\label{2.8}}
A_-:=\liminf_{n\to \infty} A(\mm_n)\ge A(\mm).
 \end{equation}
 Рассмотрим три случая:
 \begin{enumerate}
\item Пусть $A_-=\infty$, и тогда неравенство (\ref{2.8}) имеет место.

\item Пусть $|A_-|<\infty$, и тогда, не ограничивая общности, можно считать, что
 $$
   A_-=\lim_{n\to \infty}A(\mm_n),~~~\mm=\lim_{n\to \infty}\mm_n
    $$
 (если это не так, то последовательность $\mm_n$ можно заменить
 подходящей подпоследовательностью).
 Из свойства  (\ref{2.6}) и полунепрерывности снизу функции
 $A(\l,\mm)$ вытекает, что
 $$
   0\ge \liminf_{n\to \infty}A(-A(\mm_n),\mm_n)\ge A(-A_-,\mm).
 $$
 Из последнего следует $-A(\mm)\ge -A_-$, т.е. (\ref{2.8}).

\item В случае $A_-=-\infty$ можно считать, не ограничивая общности, что
  $$
 -\infty=A_-=\lim_{n\to \infty}A(\mm_n),~~~\mm=\lim_{n\to \infty}\mm_n.
    $$
 Для произвольного $N<\infty$  найдется $n_N<\infty$ такое, что
 при всех $n~\ge~n_N$ выполняется  $-A(\mm_n)\ge N$. Поскольку при любом
 фиксированном $\mm$  функция $A(\l,\mm)$ не убывает по аргументу $\l$,
 имеем в силу свойства (\ref{2.6})
 $$
   0\ge \limsup_{n\to \infty}A(-A(\mm_n),\mm_n)\ge
   \limsup_{n\to \infty}A(N,\mm_n)
$$
$$
   \ge
   \liminf_{n\to \infty}A(N,\mm_n)\ge A(N,\mm).
 $$
Таким образом, для произвольного $N < \infty$ выполнено $A(N,\mm)\leq 0$  и следовательно, $-A(\mm) \geq N .$ Последнее означает, что $A(\mm) = -\infty, $ что невозможно
 в силу свойства (\ref{2.5}), поэтому случай
$A_-=-\infty$ невозможен.
\end{enumerate}
Выпуклость и полунепрерывность снизу функции $A(\mm)$
доказаны.

Выпуклость и полунепрерывность снизу функции $A_{\gamma}(\mm)$
следуют непосредственно из определений этой функции
(\ref{lem2_1.2}) и уже установленных этих свойств для $A(\mm)$.
Утверждения \ref{th1.1i} леммы \ref{lem_th1}  доказаны.

\ref{th1.1ii}.  Выпуклость и полунепрерывность функций
$D(\aa)$, $D_{\gamma}(\aa)$ являются следствиями свойства \ref{L1} и
определения (\ref{1.3}).

Докажем компактность функции $D(\aa)$. Из определения (\ref{lem2_1.2}) и условия \ref{C} следует,
что функция $A(\mm)$
конечна и непрерывна в окрестности $\left(\0\right)_{\varepsilon}$
точки $\mm=\0$. Поэтому  найдётся $C<\infty$
такое, что для любого вектора  $\mm \in \left[ ({\bf 0})_{\varepsilon/2}\right]$ выполняется $A(\mm)\le C$.
Следовательно, отправляясь от определения функции $D(\aa)$, имеем
для любого $\aa\in \R^d$  неравенство
$$
  D(\aa)\ge  \frac{\varepsilon}{2|\aa|} \scalar{\aa}{\aa} -A\(  \frac{\varepsilon}{2|\aa|}\aa\)\ge
  \frac{\varepsilon}{2}|\aa|-C,
$$
доказывающее компактность полунепрерывной снизу функции $D(\aa)$.

Из определения (\ref{lem2_1.2}) вытекает, что функция
$A_{\gamma}(\mm)$ также ограничена в окрестности $\left[ ({\bf 0})_{\varepsilon/2}\right]$, следовательно, компактность функции $D_{\gamma}(\aa)$ устанавливается
аналогичным образом.
Утверждение \ref{th1.1ii} доказано.

\ref{th1.1iii}. Формулы (\ref{lem2_1.4}) вытекают из утверждения \ref{th1.1i}
 и свойства
 \ref{L2}.  Утверждение  \ref{th1.1iii} теоремы \ref{lem_th1}  доказано.

\ref{th1.1iv}. В силу определения  (\ref{lem2_1.2}) равенство $A_{\gamma}=A$
 эквивалентно соотношению
$$
-\gamma \le \inf_{\mm}A(\mm).
$$
 Поскольку из определения функции $D(\aa)$ следует $\inf\limits_{\mm}A(\mm)=-D(\0)$,
 то равенство $A_{\gamma}=A$
 эквивалентно неравенству $\gamma \ge D(\0)$.

\ref{th1.4ii}.  Обозначим $\widehat{D}(\aa)$ правую часть равенства (\ref{1.11}):
 \begin{equation}{\label{2.10}}
 \widehat{D}(\aa):=\inf_{\h\in [0,1]}\skk{D(\h,\aa)+\gamma(1-\h)}.
\end{equation}
 Для установления (\ref{1.11}) нам достаточно доказать тождество
 \begin{equation}{\label{2.11}}
 \widehat{D}(\aa)=A_{\gamma}^\mathfrak{Le}(\aa).
\end{equation}
 Рассмотрим функцию
$$
  X(\l,\mm):=\left\{
               \begin{array}{ll}
                 0, & \hbox{если $(\l,\mm)\in \mathcal{A}^{\le 0}$;} \\
                 \infty, & \hbox{если $(\l,\mm)\not\in \mathcal{A}^{\le 0}$.}
               \end{array}
             \right.
$$
Поскольку множество $\mathcal{A}^{\le 0}$
выпукло и замкнуто, то функция $X$ выпукла и полунепрерывна снизу, т.е.
принадлежит классу
$\mathcal{F}$. Очевидно,  что
\begin{equation}{\label{2.12}}
  D(\h,\aa)=\sup_{(\l,\mm)}\skk{\l \h+\scalar{\mm}{\aa}-X(\l,\mm)}=X^\mathfrak{Le}(\h,\aa).
\end{equation}
Определим далее функцию
$$
  V=V(\l,\mm):=\left\{
                       \begin{array}{ll}
                   0, & \hbox{если $\l\leq\gamma$;} \\
                   \infty, & \hbox{если $\l>\gamma$.}
                       \end{array}
                     \right.
$$
Легко заметить, что функция $V(\lambda,\mm)$ также выпукла и
полунепрерывна снизу, а ее преобразование Лежандра имеет вид
$$
V^\mathfrak{Le}(\theta,\aa):= \sup_{\lambda,\mm}
\skk{ \lambda\theta + \scalar{\mm}{\aa} - V(\lambda,\mm)} =\left\{
                                                \begin{array}{ll}
                   \infty, & \hbox{если $\aa\not = \0~~${\it или}$~~\h<0$;} \\
                   \gamma \h, & \hbox{если $\aa=\0~~${\it и}$~~\h\ge 0$.}
                                                \end{array}
                                              \right.
$$
Определим далее функцию
$$
  \widehat{D}(\h,\aa):=\inf_{0\le u\le \h}\skk{D(u,\aa)+\gamma(\h-u)},
  $$
и заметим (см. определение (\ref{2.10}) функции $\widehat{D}(\aa)$), что
\begin{equation}{\label{2.13}}
 \widehat{D}(\aa)=\widehat{D}(1,\aa).
\end{equation}
С другой стороны, учитывая  вид функции $V^\mathfrak{Le}(\h,\aa)$,
получаем
равенства
\begin{equation*}
\begin{split}
 \widehat{D}(\h,\aa)= &
\inf_{0\le u\le \h}\skk{D(u,\aa)+V^\mathfrak{Le}(\h-u,\0)} \\
  =& \inf_{0\le u< \infty,\bb\in \mathbb{R}^d}\skk{D(u,\bb)+V^\mathfrak{Le}(\h-u,\aa-\bb)}.
   \end{split}
\end{equation*}
Учитывая тот факт, что $D(\theta,\aa) = \infty$ для $\theta <0$ получаем
\begin{equation}\label{2.14}
 \widehat{D}(\h,\aa)=\inf_{\theta,\bb}\skk{D(u,\bb)+V^\mathfrak{Le}(\h-u,\aa-\bb)} = D * V^\mathfrak{Le}(\theta,\aa).
\end{equation}
Положим
$$
\widehat{X}(\l,\mm):=X(\l,\mm)+V(\l,\mm)=\left\{
                                           \begin{array}{ll}
 0, & \hbox{если $(\l,\mm)\in \mathcal{A}^{\le 0}_{\gamma}$;} \\
 \infty, & \hbox{если $(\l,\mm)\not\in \mathcal{A}^{\le 0}_{\gamma}$.}
                                           \end{array}
                                         \right.
$$
Тогда, используя  (\ref{2.12}), (\ref{2.14}) и свойство \ref{L3},
 получаем
$$
 \widehat{D}(\h,\aa) = X^\mathfrak{Le} * V^\mathfrak{Le}(\h,\aa)=
$$
\begin{equation}{\label{2.15}}
 (X+V)^\mathfrak{Le}(\h,\aa)=\widehat{X}^\mathfrak{Le}(\h,\aa)=
 \sup_{(\l,\mm)}\skk{\l \h+\scalar{\mm}{\aa}-\widehat{X}(\l,\mm)}.
\end{equation}

Убедимся, далее, что справедливо тождество
\begin{equation}{\label{2.16}}
  A_{\gamma}(\mm)=-\sup_{\l}\skk{\l-\widehat{X}(\l,\mm)}.
\end{equation}
Действительно,
\begin{equation*}
\begin{split}
 A_{\gamma}(\mm)= & \max\skk{-\gamma,\,A(\mm)}=\max\skk{-\gamma,\,-\sup\{\l:~(\l,\mm)\in
  \mathcal{A}^{\le 0}\}} \\
  = &-\min\skk{\gamma,\,\sup\skk{\l\,:\,(\l,\mm)\in \mathcal{A}^{\le 0}}}=
-\sup\skk{\l\,:\,(\l,\mm)\in \mathcal{A}^{\le 0}_{\gamma}}\\
= & -\sup_{\l}\skk{\l-\widehat{X}(\l,\mm)}.
\end{split}
\end{equation*}

Поэтому, учитывая (\ref{2.13}), (\ref{2.15}), (\ref{2.16}), находим
\begin{equation*}
\begin{split}
 \widehat{D}(\aa)= & \widehat{D}(1,\aa)=
 \sup_{(\l,\mm)}\skk{\l+\scalar{\mm}{\aa}-\widehat{X}(\l,\mm)} \\
  = &\sup_{\mm}\skk{\scalar{\mm}{\aa}+\sup_{\l}\skk{\l-\widehat{X}(\l,\mm)}}=
  \sup_{\mm}\skk{\scalar{\mm}{\aa}-A_{\gamma}(\mm)}.
\end{split}
\end{equation*}
Тем самым мы доказали  равенство  (\ref{2.11}), а вместе с ним и равенство
(\ref{1.11}).

\qed

\subsection{ Доказательство леммы~\ref{lem_th2}}\label{subs2.3}

\ref{th1.4i}. Равенство (\ref{1.9}) следует из определений функций $D(\aa), A(\mm)$ и  цепочки равенств
\begin{equation*}
\begin{split}
D(\aa) & =\sup_{\mm} \skk{\scalar{\mm}{\aa} - A(\mm)} =
\sup_{\mm}
\skk{\scalar{\mm}{\aa} +
\sup_{\lambda\,:\, A(\lambda,\mm) \le 0}
\{\lambda\}}\\
& =
 \sup_{(\lambda,\mm) \in \mathcal{A}^{\leq 0}}
\skk{\lambda + \scalar{\mm}{\aa}}.
\end{split}
\end{equation*}
Аналогично устанавливаются равенство в (\ref{1.10}). Пункт \ref{th1.4i} доказан.


\ref{th1.4iii}.  Для доказательства (\ref{1.13}) докажем прежде тождество
\begin{equation}{\label{2.17}}
 D_\Lambda^\mathfrak{Le}(\l,\mm)=D^\mathfrak{Le}(\l,\mm).
\end{equation}
 Имеем
\begin{equation*}
 \begin{split}
 D_\Lambda^\mathfrak{Le}(\l,\mm) & =\sup_{(\h,\aa)}\SL\l\h+\scalar{\mm}{\aa}-
 \inf_{r>0}r\Lambda\(\frac{\h}{r},\frac{\aa}{r}\)\SP\\
& =
\sup_{(\h,\aa)}\sup_{r>0}\SL\l\h+\scalar{\mm}{\aa}-
 r\Lambda\(\frac{\h}{r},\frac{\aa}{r}\)\SP \\
 & = \sup_{r>0}\skk{ r\sup_{(\h,\aa)}\skk{\l\frac{\h}{r}+  \scalar{\mm}{\frac{\aa}{r}}-
  \Lambda\(\frac{\h}{r},\frac{\aa}{r}\)}}=
  \sup_{r>0}r A(\l,\mm) \\
& =X(\l,\mm):=
      \left\{
        \begin{array}{ll}
 0, & \hbox{если $(\l,\mm)\in \mathcal{A}^{\le 0}$;} \\
 \infty, & \hbox{если $(\l,\mm)\not\in \mathcal{A}^{\le 0}$.}
        \end{array}
      \right.
\end{split}
\end{equation*}

{ Очевидно, что (см. определение $D(\theta,\aa)$  в пункте \ref{th1.4ii} теоремы \ref{lem_th1})
 $$
D(\theta,\aa) = \sup_{(\lambda,\mm)}\skk{\lambda\theta + \scalar{\mm}{\aa} - X(\lambda,\mm)}, \ \     X(\l,\mm)=D^\mathfrak{Le}(\l,\mm),
 $$
 поэтому тождество (\ref{2.17}) доказано.
 Применяя к левой и правой частям (\ref{2.17}) преобразование
 Лежандра, получаем в силу свойства \ref{L2}
 формулу (\ref{1.13}).


Докажем теперь формулу (\ref{1.14.2}). Обозначим правую часть (\ref{1.14.2}) через $\widehat{D}_{\Lambda}(\theta,\aa)$:
\begin{equation*}
  \widehat{D}_{\Lambda}(\theta,\aa) : = \lim_{\v \downarrow 0} \inf_{\aa'\in (\aa)_\v}D_{\Lambda}(\theta,\aa').
\end{equation*}
Поскольку
\begin{equation*}
  \inf_{\aa' \in (\aa)_{\varepsilon}} D_{\Lambda} (\theta, \aa') \ge \inf_{(\theta',\aa')\in ((\theta,\aa))_\v}D_{\Lambda}(\theta',\aa'),
\end{equation*}
то в силу свойства \ref{L5} имеем:
\begin{equation*}
  \widehat{D}_{\Lambda}(\theta,\aa) \ge \lim_{\varepsilon \downarrow 0 }\inf_{(\theta',\aa')\in ((\theta,\aa))_\v}D_{\Lambda}(\theta',\aa') = \left(D_{\Lambda}^{\mathfrak{Le}}\right)^{\mathfrak{Le}}(\theta,\aa).
\end{equation*}
Используя далее равенство (\ref{1.13}), получаем
\begin{equation}\label{013.16.3}
  \widehat{D}_{\Lambda}(\theta,\aa) \ge D(\theta,\aa).
\end{equation}
Далее, из определения функции $D_{\Lambda}(\theta,\aa)$ для всех $\theta', \theta >0$ имеем
\begin{equation*}
  D_{\Lambda}\left(\theta',\aa'\right) = \frac{\theta}{\theta'} D_{\Lambda}\left(\theta,\aa' \frac{\theta}{\theta'}\right),
\end{equation*}
следовательно, в силу (\ref{1.13}) и \ref{L5}
\begin{equation}\label{013.16.4}
  \begin{split}
     D(\theta,\aa) & = \lim_{\v \downarrow 0} \inf_{(\theta',\aa')\in ((\theta,\aa))_\v}D_{\Lambda}(\theta',\aa')
       = \lim_{\v \downarrow 0} \inf_{(\theta',\aa')\in ((\theta,\aa))_\v}  \frac{\theta}{\theta'} D_{\Lambda}\left(\theta,\aa' \frac{\theta}{\theta'}\right) \\
       & \ge   \lim_{\v \downarrow 0} \inf_{(\theta',\aa')\in ((\theta,\aa))_\v}  (1-\varepsilon) D_{\Lambda}\left(\theta,\aa' \frac{\theta}{\theta'}\right)
        =  \lim_{\delta \downarrow 0} \inf_{\aa'\in (\aa)_{\delta}}  D_{\Lambda}\left(\theta,\aa' \right) \\
       & =\widehat{D}_{\Lambda}(\theta,\aa).
  \end{split}
\end{equation}
Из (\ref{013.16.3}), (\ref{013.16.4}) вытекает равенство (\ref{1.14.2}).
 Пункт \ref{th1.4iii} доказан.

\ref{th1.4iv}.
Для $\h\in (0,1]$ справедливо
 $$
   D(\h,\aa)=\h D(1,\frac{\aa}{\h}).
 $$
 Докажем сначала, что для любого $\mm$
  \begin{equation}{\label{2.20}}
\widehat{D}_{\gamma}^\mathfrak{Le}(\mm)=A_{\gamma}(\mm).
\end{equation}
При $\gamma<\infty$ имеем
\begin{equation*}
\begin{split}
  \widehat{D}_{\gamma}^\mathfrak{Le}(\mm)= & \sup_{\aa}
  \skk{\scalar{\mm}{\aa}-\widehat{D}_{\gamma}(\aa)} =  \sup_{\aa}
  \skk{\scalar{\mm}{\aa}-\inf_{\h\in (0,1)}
  \skk{ \h D\(\frac{\aa}{\h}\)+\gamma(1-\h)}}  \\
= & \sup_{\aa}\sup_{\h\in (0,1)}\SL\scalar{\mm}{\aa}-
  \h D\(\frac{\aa}{\h}\)-\gamma(1-\h) \SP \\
 = & \sup_{\h\in (0,1)}\SL\h\sup_{\aa}\SL \scalar{\mm}{\frac{\aa}{\h}}-
  D\(\frac{\aa}{\h}\)\SP-\gamma(1-\h) \SP \\
  = &\sup_{\h\in (0,1)}\{\h A(\mm)-\gamma(1-\h)\}=\max\{-\gamma, A(\mm)\}=A_{\gamma}(\mm).
 \end{split}
\end{equation*}
В последнем равенстве мы использовали определение функции $A_{\gamma}(\mm)$
(см.  (\ref{lem2_1.2})).
Формула (\ref{2.20}) установлена.
Применяя к левой и правой частям (\ref{2.20}) преобразование
 Лежандра, получаем
 формулу (\ref{1.15}). Лемма~\ref{lem_th2} доказана.


\begin{thebibliography}{99}

\bibitem{MogProk3}
А.А. Могульский, Е.И. Прокопенко, Принцип больших уклонений в фазовом пространстве для многомерного первого обобщенного процесса восстановления, Сиб. электрон. матем. изв., 16, (2019),  с. 1464--1477.

\bibitem{MogProk5}
А.А. Могульский, Е.И. Прокопенко, Принцип больших уклонений в фазовом пространстве для многомерного второго обобщенного процесса восстановления, Сиб. электрон. матем. изв., 16, (2019),  с. 1478--1492.

\bibitem{Tsi}  B. Tsirelson,
From uniform renewal theorem
to uniform large and
moderate deviations for renewal-reward processes, Electron. Commun. Probab.,
18(52), (2013), pp. 1--13.

\bibitem{BorMog2}
А.А. Боровков, А.А. Могульский, Интегро--локальные предельные теоремы для обобщенных процессов восстановления при выполнении условия Крамера. I, II, Сиб. матем. журн., 59(3), (2018),  с. 491--513;
Сиб. матем. журн., 59(4), (2018),  с. 578--597.

\bibitem{MogProk1}
А.А. Могульский, Е.И. Прокопенко, Интегро--локальные теоремы для многомерных обобщенных процессов восстановления при моментном условии Крамера. I, II, III, Сиб. электрон. матем. изв., 15, (2018),     с. 475--502; Сиб. электрон. матем. изв., 15, (2018),     с. 503--527;
Сиб. электрон. матем. изв., 15, (2018),     с. 528--553.

\bibitem{MogProk2}
А.А. Могульский, Е.И. Прокопенко, Локальные теоремы для арифметических многомерных обобщенных процессов восстановления при выполнении условия Крамера, Матем. тр., 22(2), (2019),  с. 106--133.

\bibitem{MogProkLog1}
A. Logachov, A. Mogulskii, E. Prokopenko, A. Yambartsev, Local theorems for
(multidimensional) additive functionals of semi-Markov chains, Stochastic Processes and their Applications, 137, (2021), pp. 149--166.

\bibitem{MogProk4}
А.А. Могульский, Е.И. Прокопенко, Принцип больших уклонений для конечномерных распределений многомерных обобщенных процессов восстановления, Матем. тр., 23(2), (2020), с. 148--176.

\bibitem{MogLog1}
А.В. Логачёв, А.А. Могульский, Локальные теоремы для конечномерных приращений арифметических многомерных обобщенных процессов восстановления при выполнении условия Крамера, Сиб. электрон. матем. изв., 17, (2020), с. 1766--1786.


\bibitem{BorMog1}
А.А. Боровков, А.А. Могульский, Принципы больших уклонений для траектории обобщенных процессов восстановления. I, II, III, Теория вероятн. и ее примен., 60(2), (2015), с. 227--247;
Теория вероятн. и ее примен., 60(2), (2015), с. 227--247;
Теория вероятн. и ее примен., 60(3), (2015), с. 417--438.

\bibitem{MogLog2}
A.V. Logachov, A.A. Mogulskii, Anscombe–type theorem and moderate deviations for
trajectories of a compound renewal process, Journal of Mathematical Sciences, 29, (2018),  pp. 36–50.

\bibitem{Mog1}
А.А. Могульский, Расширенный принцип больших уклонений для траекторий обобщенного процесса восстановления, Матем. тр., 24(1), (2021), с. 142--174.

\bibitem{LefMarZam}  R. Lefevere,
M. Mariani, L. Zambotti,
Large deviations for renewal processes,
Stochastic
Processes and their Applications, 121(10), (2011), pp. 2243--2271.

\bibitem{Bak}
Г.А. Бакай, Большие уклонения обрывающихся многомерных обощенных
процессов востановления,
Теория вероятн. и ее примен., 66(2), (2021), с. 261--283.


\bibitem{Zamp1}
M. Zamparo,
Journal of Physics A: Mathematical and Theoretical, 52(49), (2019), 495004.

\bibitem{Zamp2}
M. Zamparo,
Large deviations in discrete-time renewal theory,
Stochastic Processes and their Applications, 139, (2021), pp. 80--109.


\bibitem{Gant}
G. Giacomin,
Random Polymer Models,  Imperial College Press, London, 2007.

\bibitem{Holl}
F. Hollander,
Random Polymers,  Springer-Verlag Berlin Heidelberg, 2009.


\bibitem{BarbuPrecupanu2012}
B. Viorel, T. Precupanu. Convexity and optimization in Banach spaces. Springer Science \& Business Media, 2012.

\bibitem{Zalinescu2002}
C. Zalinescu, Convex analysis in general vector spaces. World scientific, 2002.

\bibitem{MogPro2019}
А.А. Могульский, Е.И. Прокопенко, Функция уклонений и базовая функция для многомерного обобщенного процесса восстановления, Сиб. электрон. матем. изв., 16 (2019), с. 1449--1463.

\bibitem{Bor1}
А.А.  Боровков,  Асимптотический анализ случайных блужданий.
Быстроубывающие распределения скачков, М. "Физматлит"$,$ 2013.































\end{thebibliography}
\end{document}